\documentstyle[12pt, epsf, epsfig, amssymb, amstext, amstex]{article}
\input xy
\xyoption{all}
\begin{document}

\newtheorem{thm}{Theorem}[section]
\newtheorem{prop}[thm]{Proposition}
\newtheorem{cor}[thm]{Corollary}
\newtheorem{lem}[thm]{Lemma}
\newtheorem{conj}[thm]{Conjecture}
\newtheorem{exa}[thm]{Example}
\newtheorem{defn}[thm]{Definition}
\newtheorem{clm}[thm]{Claim}
\newtheorem{eex}[thm]{Exercise}
\newtheorem{obs}[thm]{Observation}
\newtheorem{note}[thm]{Notation}
 
\newcommand{\ben}{\begin{enumerate}}
\newcommand{\een}{\end{enumerate}}
\newcommand{\blem}{\begin{lem}}
\newcommand{\elem}{\end{lem}}
\newcommand{\bcl}{\begin{clm}}
\newcommand{\ecl}{\end{clm}}
\newcommand{\bthm}{\begin{thm}}
\newcommand{\ethm}{\end{thm}}
\newcommand{\bpr}{\begin{prop}}
\newcommand{\epr}{\end{prop}}
\newcommand{\bco}{\begin{cor}}
\newcommand{\eco}{\end{cor}}
\newcommand{\bcon}{\begin{conj}}
\newcommand{\econ}{\end{conj}}
\newcommand{\bde}{\begin{defn}}
\newcommand{\ede}{\end{defn}}
\newcommand{\bex}{\begin{exa}}
\newcommand{\eexa}{\end{exa}}
\newcommand{\bexe}{\begin{exe}}
\newcommand{\eexe}{\end{exe}}
\newcommand{\bobs}{\begin{obs}}
\newcommand{\eobs}{\end{obs}}
\newcommand{\bnote}{\begin{note}}
\newcommand{\enote}{\end{note}}

\newcommand{\fg}{\Pi _1(D-K,u)}
\newcommand{\Z}{{\Bbb Z}}
\newcommand{\C}{{\Bbb C}}
\newcommand{\R}{{\Bbb R}}
\newcommand{\Q}{{\Bbb Q}}
\newcommand{\F}{{\Bbb F}}
\newcommand{\N}{{\Bbb N}}

\newcommand{\fnref}[1]{~(\ref{#1})}
\newenvironment{emphit}{\begin{itemize} \em}{\end{itemize}}
\begin{center}
\Large{\bf {Hurwitz Equivalence in Braid Group $B_3$}}\\ 
\vspace{7mm}
\large{T. Ben-Itzhak and M. Teicher}
\footnote{Partially supported by the Emmy Noether Research Institute for Mathematics, Bar-Ilan University and the Minerva Foundation, Germany.\\
Partially supported by the Excellency Center ``Group theoretic methods in the study of algebraic varieties'' of the National Science Foundation of Israel.\\}\\
\end{center}


\vspace{15mm}
ABSTRACT. In this paper we prove certain Hurwitz equivalence properties of $B_n$. In particular we prove that for $n=3$ every two Artin's factorizations of $\Delta _3 ^2$ of the form $H_{i_1} \cdots H_{i_6}, \quad F_{j_1} \cdots F_{j_6}$ (with $i_k , j_k \in \{ 1,2 \}$) where $\{ H_1 , H_2 \} , \{ F_1 , F_2 \}$ are frames, are Hurwitz equivalent. The proof provided here is geometric, based on a newly defined frame type.\\
The results will be applied to the classification of algebraic surfaces up to deformation. It is already known that there exist surfaces that are diffeomorphic to each other but are not deformations of each other (Manetti example). We are constructing a new invariant based on Hurwitz equivalence class of factorization, to distinguish among diffeomorphic surfaces which are not deformation of each other. The precise definition of the new invariant can be found at \cite{KuTe} or \cite{Te}.
The main result of this paper will help us to compute the new invariant.

\section{Basic Definitions}
In this section we recall some basic definitions and statements from \cite{MoTe1}:\\
Let $D$ be a closed disk on $\R ^2$, $K \subset D$ finite set, $u \in \partial D$. Any diffeomorphism of $D$ which fixes $K$ and is the identity on $\partial D$ acts naturally on $\Pi _1 = \Pi _1(D-K,u)$. We say that two such diffeomorphisms of $D$ (which fix $K$ and equal identity on $\partial D$) are equivalent if they define the same automorphism on $\Pi _1(D-K,u)$. This equivalence relation is compatible with composition of diffeomorphism and thus the equivalence classes form a group.\\ 


\bde \label{BraidGroupDef}
Braid Group $B_n [D,K]$.
\ede
Let $D,K$ be as above, $n = \# K$, and let $\mathcal{B}$ be the group of all diffeomorphisms $\beta$ of $D$ such that $\beta (K)=K$, $\beta |_{\partial D} = Id_{\partial D}$. For $\beta _1, \beta _2 \in \mathcal{B}$ we say that $\beta _1$ is equivalent to $\beta _2$ if $\beta _1$ and $\beta _2$ define the same automorphism of $\Pi _1(D-K,u)$. The quotient of $\mathcal{B}$ by this equivalence relation is called the Braid group $B_n [D,K]$.\\
Equivalently, if we take the canonical homomorphism $\psi :{\mathcal{B}} \rightarrow {Aut}(\Pi _1(D-K,u))$ then $B_n [D,K] = Im(\psi )$. The elements of $B_n [D,K]$ are called braids.\\

\blem
If $K' \subset  D'$ in another pair as above with ${\# K'} = {\# K} = n$ then $B_n [D',K']$ is isomorphic to $B_n [D,K]$.
\elem

\noindent This gives rise to the definition of $B_n$:

\bde
$B_n = B_n[D,K]$ for some $D,K$ with $ \# K = n$.
\ede

\bde \label{HalfTwistDef}
$H(\sigma)$, half-twist defined by $\sigma$.
\ede
Let $D$ and $K$ be defined as above. Let $a, b \in K$, $K_{a,b} = K \backslash \{ a,b \}$ and $\sigma$ be a simple path (without a self intersection) in $D \backslash \partial D$ connecting $a$ with $b$ such that $\sigma \bigcap K = \{ a,b \}$. Choose a small regular neighborhood  $U$ of $\sigma$ such that $K_{a,b} \bigcap U = \phi$, and an orientation preserving diffeomorphism $\psi : {\R}^2 \rightarrow {\C}^1$ such that
$\psi (\sigma) = [-1,1] = \{ z \in \C ^1 | Re(z) \in [-1,1], Im(z) = 0 \}$ and $\psi (U) = \{ z \in \C ^1||z|<2 \}$. Let $\alpha (r),r \geq 0$, be a real smooth monotone function such that $\alpha (r) = 1$ for $r \in [0,3/2]$ and $\alpha (r) = 0$ for $r \geq 2$. Define a diffeomorphism $h:\C ^1 \rightarrow \C ^1$ as follows: for $z \in \C ^1, z = r e^{i\phi} \quad \text{let} \quad h(z)= r e^{i(\phi +\alpha (r) \pi )}$. It is clear that the restriction of $h$ to $\{ z \in \C ^1 | |z| \leq 3/2 \}$ coincides with the $180 ^\circ$ positive rotation, and that the restriction to $\{ z \in \C ^1 | |z| \geq 2 \}$ is the identity map. The diffeomorphism $\psi ^{-1} \circ h \circ \psi$ is inducing a braid called half-twist and denoted by $H(\sigma)$\\


\bde \label{frame}
Frame of $B_n[D,K]$
\ede

\noindent Let $K = \{ a_1,...,a_n \}$ and $\sigma _1,..., \sigma _{n-1}$ be a system of simple smooth paths in $D - \partial D$ such that $\sigma _i$ connects $a_i \text{ with } a_{i+1}$ and $L = \bigcup \sigma _i$ is a simple smooth path. The ordered system of half-twists $(H_1,...,H_{n-1})$ defined by $ \{ \sigma _i \} _{i=1} ^{n-1}$ is called a frame of $B_n [D,K]$

\bthm\label{relations} 
Let $(H_1,...,H_{n-1})$ be a frame of $B_n[D,K]$ then: 
$B_n[D,K]$ is generated by $\{ H_i \} _{i=1} ^{n-1}$ and the following is a complete list of relations:\\
$H_i H_j = H_j H_i \text{ if }|i-j|>1$ (Commutative relation)\\
$H_i H_j H_i = H_j H_i H_j \text { if } |i-j|=1$ (Triple relation)\\
\ethm

\noindent Proof can be found for example in \cite{MoTe1}.\\
This theorem provides us with Artin's algebraic definition of Braid group. \\


\section{Properties of Hurwitz Moves}



In this section we define Hurwitz Move and Hurwitz Equivalence and prove a certain theorem on the Hurwitz Equivalence of words in the Braid group which is correct for any Braid group of arbitrary order.\\


\bde
Hurwitz move on $G^m$ ($R_k, R_k ^{-1}$)
\ede
Let $G$ be a group, $\overrightarrow{t}=(t_1,...,t_m) \in G^m$. We say that $\overrightarrow{s}=(s_1,...,s_m) \in G^m$ is obtained from $\overrightarrow{t}$ by the Hurwitz move $R_k$ (or $\overrightarrow{t}$ is obtained from $\overrightarrow{s}$ by the Hurwitz move $R_k ^{-1}$) if 
$$ s_i = t_i \quad \text{ for } i \not= k, k+1,$$   
$$ s_k = t_k t_{k+1} t_k ^{-1}, \quad s_{k+1} = t_k.$$
 

\bde
Hurwitz move on factorization 
\ede
Let $G$ be a group and $t \in G$. Let $t = t_1 \cdots t_m = s_1 \cdots s_m$ be two factorized expressions of $t$. We say that $s_1 \cdots s_m$ is obtained from $t_1 \cdots t_m$ by the Hurwitz move $R_k$ if $(s_1,...,s_m)$ is obtained from $(t_1,...,t_m)$ by Hurwitz the move $R_k$.\\

\bde
Hurwitz equivalence of factorization
\ede
The factorizations $s_1 \cdots s_m$, $t_1 \cdots t_m$\label{formin2} are Hurwitz equivalent if they are obtained from each other by a finite sequence of Hurwitz moves. The notation is  $t_1 \cdots t_m \overset{HE}{\backsim} s_1 \cdots s_m$.\label{formin1} \\


\bde
Word in $B_n$
\ede
A word in $B_n$ is a representation of braid as a sequence of the frame elements and their inverses.\\\\

\noindent We will need a certain result of Garside:
\bcl\label{PositiveEqual}
{\bf (Garside):} Every two positive words (all generators with positive powers) which are equal are transformable into each other through a finite sequence of positive words, such that each word of the sequence is positive and obtained from the proceeding one by a direct application of the commutative relation or the triple relation.
\ecl
Proof: \cite{Garside}.\\

\noindent The following is obvious, never the less we give a short proof for clarification:
\bcl \label{commutes}
$G$ is a group $g_1 , g_2 \in G$:\\
1. If $g_1 g_2 = g_2 g_1$ then $g_1  \cdot g_2 \overset{HE}{\backsim} g_2 \cdot g_1$\\
2. If $g_1 g_2 g_1 = g_2 g_1 g_2$ then $g_1  \cdot g_2 \cdot g_1 \overset{HE}{\backsim} g_2 \cdot g_1 \cdot g_2$
\ecl
Proof:\\ 
1. $g_1 \cdot g_2 \overset{R_1}{\rightarrow} g_1 g_2 g_1 ^{-1} \cdot g_1 = g_2 \cdot g_1$.\\
2. $g_1 \cdot g_2 \cdot g_1 \overset{R_2}{\rightarrow} g_1 \cdot g_2 g_1 g_2 ^{-1} \cdot g_2 \overset{R_1}{\rightarrow} g_1 g_2 g_1 g_2 ^{-1} g_1 ^{-1} \cdot g_1 \cdot g_2$ but $g_1 g_2 g_1 = g_2 g_1 g_2$ so $g_1 g_2 g_1 g_2 ^{-1} g_1 ^{-1} = g_2$ and we get $g_1  \cdot g_2 \cdot g_1 \overset{HE}{\backsim} g_2 \cdot g_1 \cdot g_2$\\

\noindent From Garside and the above we get the following proposition:

\bpr \label{sequence}
Let $H_1,...,H_{n-1}$ be a set of generators of $B_n$ and $H_{i_1} \dots H_{i_p} = H_{j_1} \dots H_{j_p}$ two positive words (with $i_k, j_k \in \{ 1,...n-1 \} $) then  $H_{i_1} \cdots H_{i_p} \overset{HE}{\backsim} H_{j_1} \cdots H_{j_p}$.
\epr

\noindent Proof:\\
Apply \ref{PositiveEqual} on $H_{i_1} \dots H_{i_p} = H_{j_1} \dots H_{j_p}$ we get a finite sequence of positive words $\{ W_r \} _{r=0} ^{q}$ s.t. $W_0 = H_{i_1} \dots H_{i_p}$, $W_q = H_{j_1} \dots H_{j_p}$ and $W_{r+1}$ is obtained from $W_r$ by a single application of the commutative relation or the triple relation.\\
An application of the triple relation (as in \ref{commutes}) is equal to an application of 2 Hurwitz moves $R_{t+1}, R_t$. An application of the commutative relation is equal to an application of one Hurwitz move (on the commuting elements).
Thus, $ W_0 \overset{HE}{\backsim} W_q$

\bde
$\Delta _n ^2  \in B_n[D,K]$
\ede
$\Delta _n ^{2} = (H_1 \dots H_{n-1})^n$ where $ \{ H_i \} _{i=1} ^{n-1}$ is a frame.\\

\noindent We apply \ref{sequence} on $\Delta _n ^2$: 

\bco
All $\Delta _n ^2$ factorizations $H_{i_1} \cdots H_{i_{n(n-1)}} \quad i_k \in \{ 1,...,n-1\} $ are Hurwitz Equivalent.
\eco


\section{$B_3$ Properties}

This section concentrates on $B_3$ and we prove the following new result:

\bthm
Let $\{ H_1, H_2 \}, \{ F_1, F_2 \}$ two frames of $B_3$, then for every two words $\Delta ^2 _3 = H_{i_1} ... H_{i_6} = F_{j_1} ... F_{j_6}$ (with $i_k,j_k \in \{ 1,2 \}$) the factorizations $H_{i_1} \cdots H_{i_6}$,  $F_{j_1} \cdots F_{j_6}$ are Hurwitz equivalent.
\ethm
The proof of the main theorem is at the end of the section following from propositions \ref{partial} and \ref{reduce}.\\\\

\begin{figure}[h]
\begin{center}
\epsfxsize=4cm
\epsfysize=4cm
\epsfbox{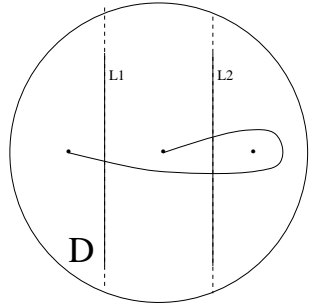}
\caption{Simple path with $m(\sigma ) = 3$}
\label{sigLength}
\end{center}
\end{figure}

\noindent Let $D$ be a closed disk on $\R ^2$ and $K \subset D$, $K=\{ a,b,c \}$. We fix two parallel lines $L_1, L_2$ in $D$  which separate the points in $K$ (see figure \ref{sigLength}).

\noindent With the choice of $L_1$ and $L_2$ we define the length and the minimal paths related to a half-twist:
\bde
Let $\sigma $ be a simple path in $D \backslash \partial D$ that defines the half-twist $H(\sigma )$. Then $m(\sigma )=  \# \{ \sigma \bigcap (L_1 \bigcup L_2) \} $ as demonstrated in figure \ref{sigLength}.
\ede

\bde
Let $H$ be a half-twist in $B_n[D,K]$. We define the length of $H$, $m(H) = min \{ m(\sigma ) | \quad H = H(\sigma) \}$ and we call such a path $\sigma$ which gives the minimal value, the minimal path of $H$.
\ede

\bnote \label{minFrame}
Any frame $ \{ H_1, H_2 \} $ can be defined by the minimal paths $\sigma _1, \sigma _2$ which also satisfy the relations in \ref{frame}. 
\enote

\begin{figure}[htp]
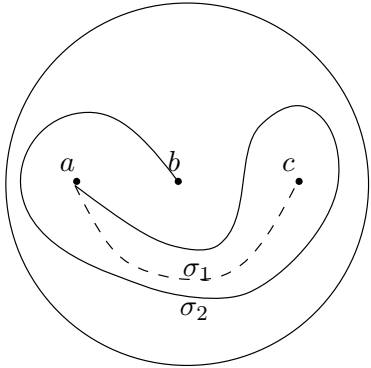

\begin{center}
\epsfxsize=8cm
\epsfysize=4cm
\input partial.pstex_t
\caption[]{
\label{partialHT}
\sf
Economic frame $m(H(\sigma _1)) < m(H(\sigma _2))$.
}
\end{center}
\end{figure}

\bde
We say that a frame $\{ H_1 , H_2 \}$ is economic if there exist $\sigma _1,\sigma _2$ minimal paths of $H_1,H_2$ respectively s.t. for every small neighborhood $U$ of $(\sigma _1\bigcap K) \backslash \sigma _2$, $\sigma _1 \subset \sigma _2$ outside $U$ (see \ref{partialHT})
\ede

\blem \label{neq}
Let $\{ H_1, H_2\} $ be economic frame then $m(H_1) \neq m(H_2)$
\elem
\noindent Proof: Trivial.

\bpr \label{partial}
Any frame $\{ H_1, H_2 \} $ which generates $B_3[D,K]$ s.t $max(m(H_1),m(H_2)) > 1$ is economic.
\epr

\noindent Proof:\\
Suppose that there exists a frame $\{ H_1, H_2 \} $ as above which is not economic. From notation \ref{minFrame}, there exist $\sigma _1 , \sigma _2$ as defined in \ref{frame} s.t. $\sigma _1 , \sigma _2$ are minimal paths of $H_1, H_2$ respectively and the two paths are 'separated' at some point.\\ 
We will show that such a 'separation' is not possible. Starting from ${s} = \sigma _1 \bigcap \sigma _2$ we will examine the first 'separation' where the two paths are separated by one of the points in $K$.\\
The cases where $s = a$ and $s = c$ are symmetric. Proving the case where $s = b$ is similar.\\
We will assume that both paths start at $a$ 

\begin{figure}[h]
\begin{center}
\epsfxsize=4cm
\epsfysize=4cm
\epsfbox{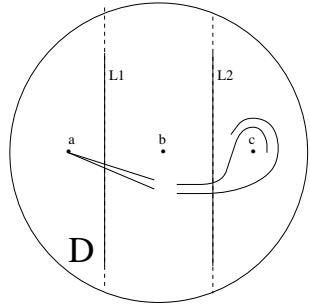}
\caption{Same point different orientation and turning back.}
\label{split1A}
\end{center}
\end{figure}

\begin{figure}[h]
\begin{center}
\epsfxsize=4cm
\epsfysize=4cm
\epsfbox{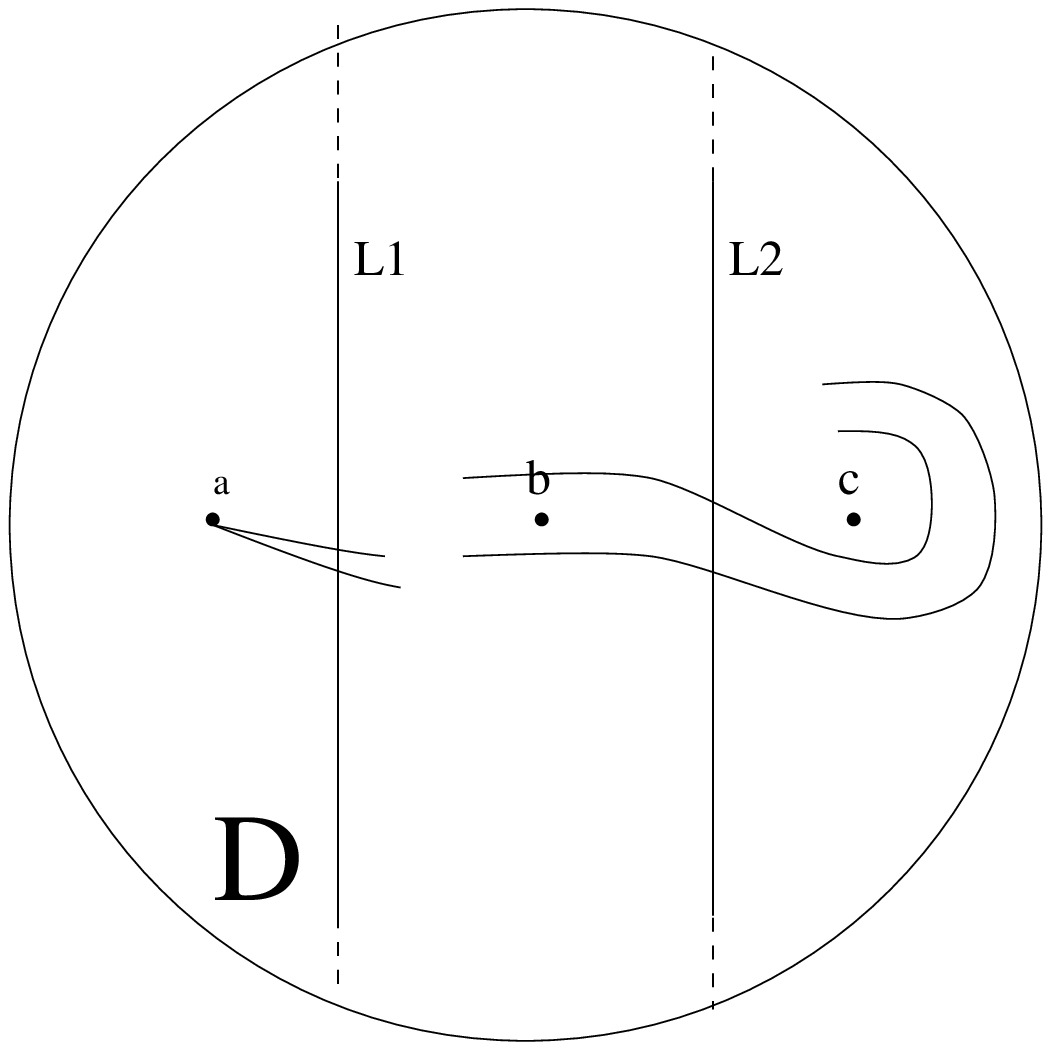}
\caption{Same point different orientation and continue to next point with the same orientation.}
\label{split1B1}
\end{center}
\end{figure}

\begin{figure}[h]
\begin{center}
\epsfxsize=4cm
\epsfysize=4cm
\epsfbox{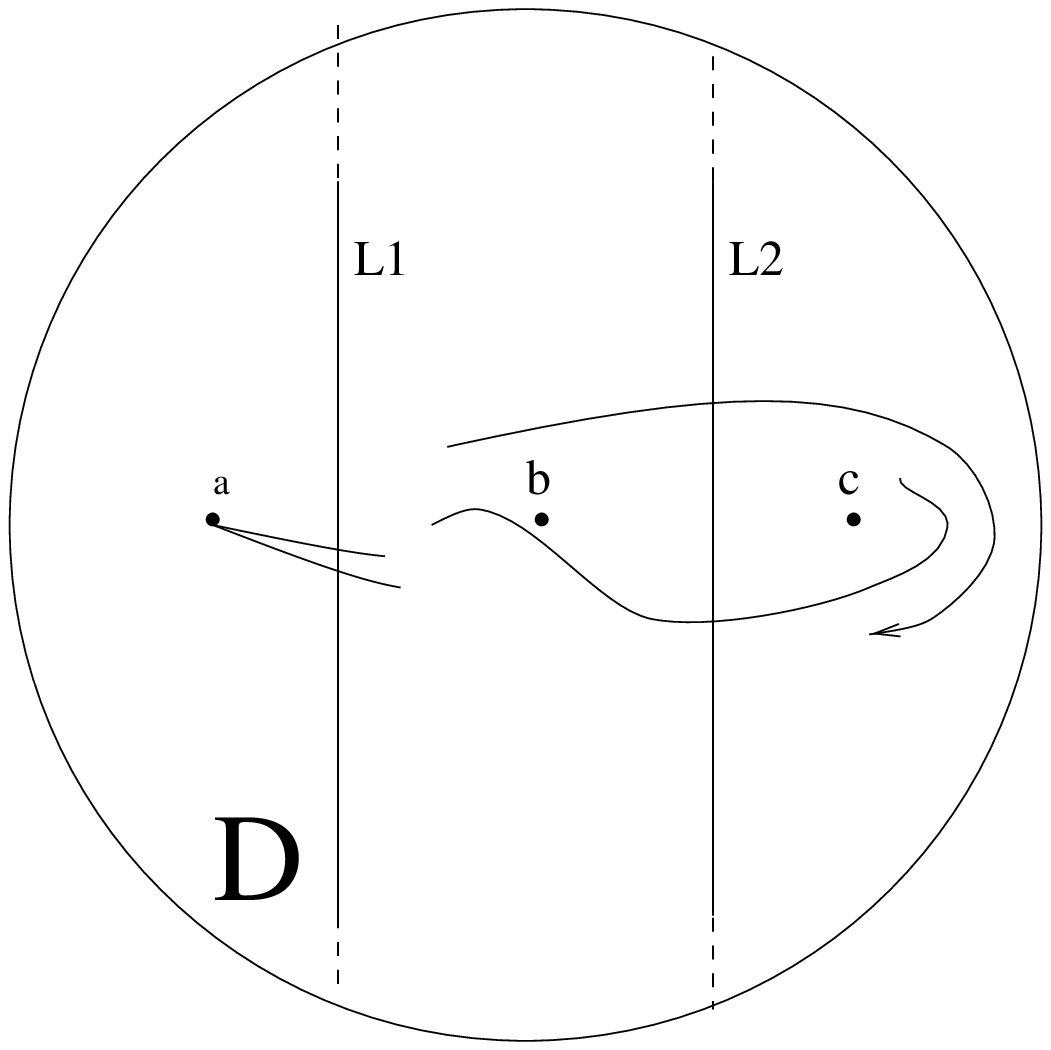}
\caption{Same point different orientation and continue to next point with different orientation.}
\label{split1B2}
\end{center}
\end{figure}

\noindent {\bf Case 1:} The two paths entering a neighborhood of the same point but with different orientation (one clockwise and the other counter-clockwise).\\

\noindent (1.a) Entering and turning back as shown in figure \ref{split1A}. The two paths create a cycle with one point inside, the internal path has no where to turn, and thus cannot end. This case is also relevant when one path is turning back and the other continues to the next point.\\

\noindent (1.b) Entering and continue to the next point:\\
(1.b.i)  The two paths continue to the next point with the same orientation (figure \ref{split1B1}). The two paths create a cycle that starts at $a$ and contains $b$. One path has to end at $b$ and will have to change it's orientation through $c$. This will close a cycle with all the points inside, and the external path has no where to end.\\
(1.b.ii) The two paths continue to the next point with different orientation (figure \ref{split1B2}). The two paths create a cycle with all points inside so the external path has no where to end.\\

\begin{figure}[h]
\begin{center}
\epsfxsize=4cm
\epsfysize=4cm
\epsfbox{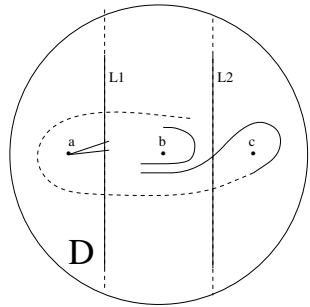}
\caption{External path must create a cycle with the internal.}
\label{split2A}
\end{center}
\end{figure}

\noindent {\bf Case 2:} Entering to the same point with the same orientation and leaving to different points.\\

\noindent This can only happen in point $b$. The internal path must go back to $a$ and the external will go to $c$.\\
\noindent (2.a) The external path goes to $c$ by changing orientation (figure \ref{split2A}). This enforces the external path to go around $a$ and create a cycle with the internal path. Thus, the two paths create a cycle with all points inside and the external path has no where to end.\\

\begin{figure}[h]
\begin{center}
\epsfxsize=4cm
\epsfysize=4cm
\epsfbox{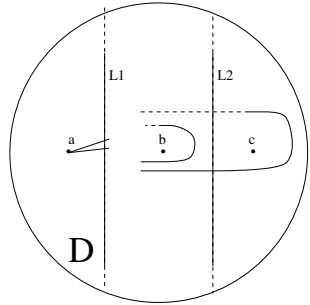}
\caption{External path must create a cycle with the internal.}
\label{split2B}
\end{center}
\end{figure}

\noindent (2.b) The external path goes to $c$ with the same orientation (figure \ref{split2B}). The external path cannot go to the common part (this will create an empty cycle and the path has no where to end) of the paths so both paths continue to $a$. In this case we get an equivalent cycle as in 1.b.i.\\
This conclude the proof of proposition \ref{partial}.

\blem \label{conjFrame}
Let $\{ H_1, H_2 \} $ be a frame generates $B_3 [D,K]$ and $H_\zeta$ half-twist, than  $\{ H {_\zeta} ^{-1} H_1 H_{\zeta},  H {_\zeta} ^{-1} H_2 H_{\zeta} \} $ is also a frame.
\elem
\noindent Proof:\\
By III.1.0 of \cite{MoTe5},
$$H_1 ^{-1} H_2 H_1 = H(H_1(\sigma _2 )) \text{ where } H_2 = H(\sigma _2).$$
$\{ H_1 = H(\sigma _1), H_2 = H(\sigma _2) \} $ is a frame, therefore $\sigma _1 \bigcup \sigma _2$
is a simple smooth path. $H_{\zeta} $ is a diffeomorphism so $H_{\zeta} (\sigma _1) \bigcup H_{\zeta}(\sigma _2) = H_{\zeta} (\sigma _1 \bigcup \sigma _2)$ is also a simple smooth path.\\
$K \subset \sigma _1 \bigcup \sigma _2, \quad H_{\zeta}(K)=K$ and therefore $H_{\zeta} (\sigma _1), H_{\zeta}(\sigma _2)$ connecting all points in $K$ appropriately.

\bpr \label{reduce}
Let $\{ H_1, H_2 \} $ be a frame of $B_3 [D,K]$ then there exist a frame $\{ F_1, F_2 \}$ s.t. 
$$H_1 \cdot H_2 \cdot H_1 \cdot H_2 \cdot H_1 \cdot H_2 \overset{HE}{\backsim} F_1 \cdot F_2 \cdot F_1 \cdot F_2 \cdot F_1 \cdot F_2$$
$$ \text{and } max(m(H_1), m(H_2)) > max (m(F_1), m(F_2)).$$
\epr

\noindent Proof:\\
From proposition \ref{partial} we get that $\{ H_1,H_2\} $ is economic and from lemma \ref{neq} $m(H_1) \neq m(H_2)$ we will assume that $m(H_1)<m(H_2)$.\\

\begin{figure}[htp]
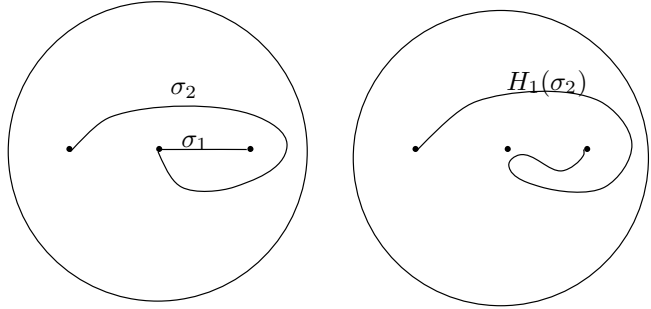

\begin{center}
\epsfxsize=8cm
\epsfysize=4cm
\input conj.pstex_t
\caption[]{
\label{conj}
\sf
Half-twist conjugation.
}
\end{center}
\end{figure}

\noindent Since $\{ H_1,H_2\} $ is economic and $m(H_1)<m(H_2)$, operating $H_1$ on $\sigma _2$ (conjugation) will cause $\sigma _2$ to start coming back from the starting point ($\sigma _1 \bigcap \sigma _2 \bigcap K$) along the common part ($\sigma _1$) as shown in figure \ref{conj}. Therefore, we get that either
$$m(H _1 ^{-1} H_2 H_1) = m(H(H_1(\sigma _2 ))) = m(H_2) -m(H_1) \text{        or }$$ 
$$m(H _1 H_2 H_1 ^{-1}) = m(H(H_1 ^{-1}(\sigma _2 ))) = m(H_2) -m(H_1)$$

\noindent In the case where $m(H(H_1(\sigma _2 ))) = m(H_2) -m(H_1)$:\\
By performing the Hurwitz moves $R_1 R_0 R_4 R_3$ we get:
$$H_1 \cdot H_2 \cdot H_1 \cdot H_2 \cdot H_1 \cdot H_2 \overset{HE}{\backsim}  H_2 \cdot H_1 \cdot H_2 \cdot H_1 \cdot H_2 \cdot H_1$$
\noindent and by performing $R_0 ^{-1} R_2 ^{-1} R_4 ^{-1}$:
$$H_2 \cdot H_1 \cdot H_2 \cdot H_1 \cdot H_2 \cdot H_1 \overset{HE}{\backsim} H_1 \cdot H_1 ^{-1} H_2 H_1 \cdot H_1 \cdot H_1 ^{-1} H_2 H_1 \cdot H_1 \cdot H_1 ^{-1} H_2 H_1$$

\noindent Since $ m(H_1 ^{-1} H_2 H_1) = m(H(H_1(\sigma _2 ))) $ and $ m(H_2) > m(H_1)$,
$$max(m(H_1), m(H_2)) > max(m(H_1), m(H_1 ^{-1} H_2 H_1))$$
From Lemma \ref{conjFrame}, the set $\{ H_1,H_1 ^{-1} H_2 H_1 \}$ is a frame, we get an equivalent factorization of the same form but with a reduced length.\\\\
In the case where $m(H(H_1 ^{-1}(\sigma _2 ))) = m(H_2) -m(H_1)$:\\
$$H_1 \cdot H_2 \cdot H_1 \cdot H_2 \cdot H_1 \cdot H_2 \overset{HE}{\backsim} H_1 H_2 H_1 ^{-1} \cdot H_1 \cdot H_1 H_2 H_1 ^{-1} \cdot H_1 \cdot H_1 H_2 H_1 ^{-1} \cdot H_1$$ By performing the Hurwitz moves $R_0 R_2 R_4$. Once again, we get a frame where the maximum length is reduced.\\\\

\noindent Proof of the main theorem:\\
Let  $H_{i_1} \cdot H_{i_2} \cdot H_{i_3} \cdot H_{i_4} \cdot H_{i_5} \cdot H_{i_6}$ and $F_{j_1} \cdot F_{j_2} \cdot F_{j_3} \cdot F_{j_4} \cdot F_{j_5} \cdot F_{j_6}$ be two factorizations as above.\\
From \ref{sequence} we get that:\\
$$H_{i_1} \cdot H_{i_2} \cdot H_{i_3} \cdot H_{i_4} \cdot H_{i_5} \cdot H_{i_6} \overset{HE}{\backsim} H_1 \cdot H_2 \cdot H_1 \cdot H_2 \cdot H_1 \cdot H_2 \quad\ \text{and} $$
$$F_{j_1} \cdot F_{j_2} \cdot F_{j_3} \cdot F_{j_4} \cdot F_{j_5} \cdot F_{j_6} \overset{HE}{\backsim} F_1 \cdot F_2 \cdot F_1 \cdot F_2 \cdot F_1 \cdot F_2$$

\noindent By \ref{reduce},both $H_1 \cdot H_2 \cdot H_1 \cdot H_2 \cdot H_1 \cdot H_2$ and $F_1 \cdot F_2 \cdot F_1 \cdot F_2 \cdot F_1 \cdot F_2$ are Hurwitz equivalent to $X_1 \cdot X_2 \cdot X_1 \cdot X_2 \cdot X_1 \cdot X_2$ where $\{ X_1 , X_2 \}$ is a frame with $m(X_1 ) = m(X_2) =1$. Since such a frame is unique we get that 
$$H_{i_1} \cdot H_{i_2} \cdot H_{i_3} \cdot H_{i_4} \cdot H_{i_5} \cdot H_{i_6} \overset{HE}{\backsim} F_{j_1} \cdot F_{j_2} \cdot F_{j_3} \cdot F_{j_4} \cdot F_{j_5} \cdot F_{j_6}$$

\end{document}